\documentclass[12pt]{article}
\usepackage{amssymb,latexsym,a4}
\newcommand{\norm}[2]{{\mathrm{N}}_{#1}\left({#2}\right)}
\newcommand{\into}{\hookrightarrow}
\newcommand{\End}{\mathop{\mathrm{End}}\nolimits}
\newcommand{\Q}{\ensuremath{\mathbb{Q}}}
\newcommand{\GL}{\mathop{\mathrm{GL}}\nolimits}
\newcommand{\ebar}{\overline{e}}
\newcommand{\Z}{\ensuremath{\mathbb{Z}}}
\newcommand{\kbar}{\ensuremath{\overline{k}}}
\newtheorem{Theorem}{Theorem}
\newtheorem{Lemma}{Lemma}
\newcommand{\eps}{\varepsilon}

\newcommand{\im}{\mathop{\mathrm{im}}\nolimits}
\newenvironment{Proof}{\noindent\textbf{Proof:}}{\hspace*{\fill}$\Box$}
\newenvironment{EnumProof}{\noindent\textbf{Proof:}\begin{enumerate}}%
  {\hspace*{\fill}$\Box$\end{enumerate}}
\newcommand{\mytab}{\hspace{-3pt} & \hspace{-3pt}}
\newcommand{\F}[1]{\ensuremath{\mathbb{F}_{#1}}}
\newcommand{\mO}{\ensuremath{\mathcal{O}}}

\title{Degrees of polarizations on an abelian surface with real multiplication}
\newcommand{\MSC}{\setcounter{footnote}{1}\footnotetext{\textit{Mathematics
Subject Classification} (2000): 11G10, (11G15, 14K22)}}
\author{John Wilson}
\date{16th March 2000}
\begin{document}
\maketitle\MSC

\begin{abstract}
Let \(F\) be a real quadratic field, and let \(R\) be an order in
\(F\). Suppose given a polarized abelian surface \((A,\lambda)\)
defined over a number field \(k\) with a symmetric action of \(R\)
defined over \(k\). This paper considers varying \(A\) within the
\(k\)-isogeny class of \(A\) to reduce the degree of \(\lambda\) and
the conductor of \(R\). It is proved, in particular, that there is a
\(k\)-isogenous principally polarized abelian surface with an action
of the full ring of integers of \(F\) when \(F\) has class number
\(1\) and the degree of \(\lambda\) and the conductor of \(R\) are odd
and coprime.
\end{abstract}

\section{Introduction}

We shall say that an abelian variety \(A\) has \emph{real
multiplication} by a totally real field \(F\) of degree \(\dim A\) if
there is an embedding \(F\into\End(A)\otimes\Q\). We shall say that
\(A\) has \emph{maximal} real multiplication by \(F\) if the embedding
\(F\into\End(A)\otimes\Q\) induces an embedding of the ring of
integers of \(F\) into the endomorphism ring of \(A\).

Let \(J_0(N)\) denote the Jacobian of the modular curve \(X_0(N)\),
which para\-metr\-izes pairs of \(N\)-isogenous elliptic curves; then
the simple factors of \(J_0(N)\) always have real
multiplication. Indeed, abelian varieties with real multiplication are
examples of abelian varieties of \(\GL_2\)-type, and it is conjectured
that any abelian variety of \(\GL_2\)-type is a simple factor of
\(J_0(N)\) for suitable \(N\) \cite{RibetKorea}.

Now the simple factors of \(J_0(N)\) are determined up to isogeny, and
so the question naturally arises as to what extent one may be able to
reduce the degree of a given polarization on such a variety, and to
what extent one may be able to enlarge the endomorphism ring, within
this isogeny class.

Taylor and Shepherd-Barron have proved in \cite{T--SB} that many
abelian surfaces with real multiplication by \(\Q(\sqrt{5})\) defined
over \Q\ are indeed factors of \(J_0(N)\) for some \(N\): their
theorem requires a few technical conditions, but also requires that
the surface admits a principal polarization, and maximal real
multiplication over \Q. However, with reference to Wang's table in
\cite{Wang}, the natural polarization on a simple two-dimensional
factor of \(J_0(N)\) is rarely principal, and these factors do not
always have maximal real multiplication.

In general, principally polarized abelian surfaces are known to be
either products of elliptic curves, or Jacobians of curves of genus
2. Moreover, if \(C\) is a curve of genus 2 (defined over some number
field), then the condition that the Jacobian \(J\) of \(C\) should
have a symmetric action of the full ring of integers of a specified
field \(F\) can be characterized by the existence of certain curves on
the associated Kummer surface (the quotient of \(J\) by the involution
\(-1\)). This latter point is essentially contained in work of Humbert
\cite{Humbert} at the end of the nineteenth century; see also chapter
XVIII of Hudson's book \cite{Hudson} and, for an account in more
modern language, chapter IX of \cite{vdG}.

In view of all this, one might hope that any abelian surface \(A\)
with real multiplication by a field \(F\) is isogenous to one with a
principal polarization, or that \(A\) is isogenous to a surface with
maximal real multiplication by \(F\). We aim for a result along these
lines, and prove the following.

\begin{Theorem}\label{thm:main}
Let \(R\) be an order in \(F\) and let \((A,\lambda)\) be a polarized
abelian surface over \(k\) with a symmetric action of \(R\) defined
over \(k\). Suppose that \(F\) has class number \(1\), and that
the degree of \(\lambda\) and the conductor of \(R\) are odd and
coprime.

Then there is a \(k\)-isogenous abelian surface \(B\) with a principal
polarization \(\lambda'\) on \(B\) and a symmetric action of the ring
of integers of \(F\), both defined over \(k\).
\end{Theorem}

Our technique actually deals with increasing the endomorphism ring and
decreasing the degree of the polarization as two separate issues. The
results we obtain in each case are theorem \ref{thm:max} (section
\ref{sec:max}) and theorem \ref{thm:sr} (section \ref{sec:sr})
respectively.

\bigskip\noindent\textbf{Acknowledgement:}
I should like to thank Prof Birch for his supervision of the doctoral
thesis \cite{Thesis} of which this work forms a part.

\section{Lemmata}

Let us fix some notation. Throughout the following, \(k\) will denote
a number field, \(\kbar\) an algebraic closure of \(k\) and \(F\) a
real quadratic field.

Let \(A\) be an abelian variety defined over \(k\). We shall denote
the Weil pairings by
\begin{eqnarray*}
\ebar_m\mytab\colon\mytab A[m]\times A^{\vee}[m]\to \mu_m \\
\textrm{and }
e_\ell\mytab\colon\mytab T_\ell A\times T_\ell A^{\vee}\to\Z_\ell(1).
\end{eqnarray*}
When \(\lambda\colon A\to A^{\vee}\) is a homomorphism (in particular,
when \(\lambda\) is a polarization on \(A\)), we define two more
pairings by
\begin{eqnarray*}
\ebar_m^\lambda = \ebar_m\circ(1\times\lambda)\mytab\colon\mytab
A[m]\times A[m]\to\mu_m \\
\textrm{and }
e_\ell^\lambda = e_\ell\circ(1\times\lambda)\mytab\colon\mytab
T_\ell A\times T_\ell A\to\Z_\ell(1).
\end{eqnarray*}
See \cite[\S 16]{MilneAb} for definitions and properties of these
pairings. (Note that we shall make frequent reference to
Milne's article \cite{MilneAb} for standard facts and definitions
concerning abelian varieties.)

If \((A,\lambda)\) is a polarized abelian variety over \(k\), we
define a pairing \(e^\lambda\) on \(\ker\lambda\) as follows:
given \(a,a'\in \ker\lambda\), choose \(m\) such that \(ma=0=ma'\) and
\(b\in A(\kbar)\) such that \(mb=a'\), then set \(e^\lambda(a,a') =
\ebar_m (a,\lambda(b))\). This gives a well-defined alternating
pairing. (Again, see \cite[\S 16]{MilneAb}.)

We shall denote by \(\End_k(A)\) the ring of endomorphisms of \(A\)
which are defined over \(k\). For a polarized abelian variety
\((A,\lambda)\), we shall denote by \(\End_k(A)^\mathrm{s}\) the ring
of symmetric endomorphisms over \(k\), that is,
\[
\End_k(A)^\mathrm{s} =
\{\alpha\in\End_k(A) \mid \lambda\circ\alpha=\alpha^{\vee}\circ\lambda\}
\]

Suppose that \(\pi\colon A\to B\) is an isogeny of abelian varieties,
and that \(\lambda\) is a polarization on \(B\). Then \(\lambda\)
induces, via \(\pi\), a polarization \(\pi^{*}\lambda\) on \(A\) as
follows: the isogeny \(\pi\) induces an isogeny \(\pi^{\vee}\colon
B^{\vee}\to A^{\vee}\) of the dual abelian varieties by functoriality;
the polarization \(\pi^{*}\lambda\) is the composition
\(\pi^{\vee}\circ\lambda\circ\pi\). (Recall that the polarization
\(\lambda\) is given by the algebraic equivalence class of
some ample divisor \(D\) on \(B\). Then \(\pi^{*}\lambda\) is the
polarization on \(A\) corresponding to the preimage \(\pi^{*}D\) on
\(A\).)

If \((A,\lambda)\) and \((B,\lambda')\) are polarized abelian
varieties, then by an \emph{isogeny of polarized abelian varieties}
\(\pi\colon(A,\lambda)\to(B,\lambda')\), we shall mean an isogeny
\(\pi\colon A\to B\) such that \(\pi^{*}\lambda' = \lambda\).

\bigskip
Now we prepare the ground with a few results, some of which are well
known.

\begin{Lemma}
\label{lem:easy}
Let \(V\) be a vector space over a field \(K\) and suppose that the
characteristic of \(K\) is not \(2\). Suppose given a nondegenerate
antisymmetric and \(K\)-bilinear pairing
\(\langle{\cdot},{\cdot}\rangle\colon V\times V\to K\), and a nonzero
endomorphism \(\eps\) of \(V\) such that \(\langle u,\eps v\rangle =
\langle \eps u, v \rangle\) for all \(u,v\in V\). Then
\begin{enumerate}
\item[\textup{(i)}] \(\eps(v)\in v^\bot\) for all \(v\in V\) and
\item[\textup{(ii)}] \(\dim(\im\eps)\) is even.
\end{enumerate}
\end{Lemma}

\begin{EnumProof}
\item[(i)] Given \(v\in V\), we have \(\langle\eps(v),v\rangle =
\langle v,\eps(v)\rangle = -\langle\eps(v),v\rangle\).
\item[(ii)] We may define a pairing
\(({\cdot},{\cdot})\colon\im\eps\times\im\eps\to V\) as follows:
given \(u,v\in\im\eps\), choose \(w\in\eps^{-1}(u)\) and set
\((u,v)=\langle w,v\rangle\). One readily checks that this pairing is
well-defined, \(K\)-bilinear, antisymmetric and nondegenerate.
\end{EnumProof}

\begin{Lemma}\label{lem:desc}
Let \(\pi\colon(A,\lambda)\to(B,\lambda')\) be an isogeny of polarized
abelian varieties defined over \(k\), and suppose that
\(R\subseteq\End_k(A)^\mathrm{s}\) preserves \(\ker\pi\). Then
\(\pi\) induces \(R\into\End_k(B)^\mathrm{s}\).
\end{Lemma}

\begin{Proof}
Let \(\alpha\) be an element of \(R\). Since \(R\) preserves
\(\ker\pi\), we may choose \(\alpha'\in\End(B)\) such that
\(\alpha'\circ\pi = \pi\circ\alpha\). It is now straightforward to
check that \(\alpha'\) is symmetric and defined over \(k\) (since
\(k\) is perfect), and that the map \(\alpha\mapsto\alpha'\) is an
embedding.
\end{Proof}

\begin{Lemma}\label{lem:pol}
Let \((A,\lambda)\) be a polarized abelian variety over \(k\).

\begin{enumerate}
\item[\textup{(i)}] A homomorphism \(\lambda\colon A\to A^{\vee}\) is a
polarization if and only if the pairings \(e_\ell^\lambda\) are
antisymmetric.
\item[\textup{(ii)}] Suppose that \(\alpha\in\End_k(A)^\mathrm{s}\)
is an isogeny. Then \(\lambda = \lambda'\alpha\) for some polarization
\(\lambda'\) on \(A\) defined over \(k\) if and only if
\(\ker\lambda\supseteq\ker\alpha\).
\item[\textup{(iii)}] Let \(\pi\colon A\to B\) be an isogeny defined
over \(k\). Then \(\lambda = \pi^{*}\lambda'\) for some polarization
\(\lambda'\) on \(B\) defined over \(k\) if and only if
\(\ker\lambda\supseteq\ker\pi\) and the pairing \(e^\lambda\) is trivial on
\(\ker\pi\times\ker\pi\).
\end{enumerate}
\end{Lemma}

\begin{EnumProof}
\item[(i)] Milne proves this over algebraically closed fields of
characteristic other than \(2\) in Proposition 16.6 of
\cite{MilneAb}. In view of our hypothesis that \(k\) is a number
field, and hence perfect, Remark 16.14, \textit{loc. cit.}, applies
and the result holds in our case.
\item[(ii)] Necessity is clear.

For sufficiency, note that we can at least define a homomorphism
\(\lambda'\colon A\to A^{\vee}\), defined over \(k\), such that
\(\lambda=\lambda'\alpha\).

Now let \(a=(a_n)\in T_\ell A\). By the definitions,
\[
\ebar_{\ell^n}^{\lambda'}(a_n,a_n) = \ebar_{\ell^n}(a_n,\lambda'(a_n))
= \ebar_{\ell^n}(a_n,\lambda(b_n)),
\]
where \(b=(b_n)\in T_\ell A\) satisfies \(\alpha(b_n) = a_n\).

There is an isogeny \(\alpha'\in\End_k(A)^\mathrm{s}\) such that
\(\alpha'\alpha = \deg\alpha\). Suppose that \(\deg\alpha =
\ell^r m\), where \(\ell\) does not divide \(m\). Then
\[
\ebar_{\ell^n}(a_n,\lambda(b_n)) = 
  \ebar_{\ell^{n+r}}(a_{n+r},\lambda(b_{n+r}))^{\ell^r},
\]
and so
\begin{eqnarray*}
\ebar_{\ell^n}^{\lambda'}(a_n,a_n)^m \mytab = \mytab
  \ebar_{\ell^{n+r}}\big(a_{n+r},\lambda(b_{n+r})\big)^{(\deg\alpha)} \\
\mytab = \mytab
  \ebar_{\ell^{n+r}}(a_{n+r},\lambda(\alpha'\alpha b_{n+r})) \\
\mytab = \mytab
  \ebar^\lambda_{\ell^{n+r}}(a_{n+r},\alpha'(a_{n+r})) \\
\mytab = \mytab 1,
\end{eqnarray*}
where the last step applies lemma \ref{lem:easy}(i).

Applying part (i) of this lemma, \(\lambda'\) is a polarization.
\item[(iii)] \cite[proposition 16.8]{MilneAb}.
\end{EnumProof}

\bigskip\noindent
Note that in lemma \ref{lem:pol}(ii), we have
\(\deg\lambda = \deg\alpha\cdot\deg\lambda'\), and in (iii), we have
\(\deg\lambda = (\deg\pi)^2 \deg\lambda'\). 

\bigskip
The next lemma considers reducing the degree of a polarization by
taking isogenies without reference to any real multiplication.

\begin{Lemma}\label{lem:sqf}
Let \((A,\lambda)\) be a polarized abelian variety defined over \(k\).
Then there is a \(k\)-isogenous abelian variety \(B\) with a
polarization \(\lambda'\) on \(B\) such that:
  \begin{enumerate}
  \item[\textup{(i)}] \(\deg\lambda'\) divides \(\deg\lambda\);
  \item[\textup{(ii)}] for each prime \(\ell\) which divides
  \(\deg\lambda'\), the subgroup \((\ker\lambda')[\ell^\infty]\) is
  properly contained in \(B[\ell]\) and is of even dimension over
  \(\F{\ell}\); and
  \item[\textup{(iii)}]
  \(\End_k(A)^\mathrm{s}\into\End_k(B)^\mathrm{s}\) via \(\pi\).
  \end{enumerate}
\end{Lemma}

\begin{Proof}
Write \(\Lambda=\ker\lambda\). We shall induct on \(\deg\lambda\),
taking as our base cases those where
\(\Lambda[\ell^\infty]\) is properly contained in \(A[\ell]\) for each
\(\ell\).

In general, suppose that \(\ell\) is a prime dividing \(\deg\lambda\). If
\(\Lambda[\ell^\infty]\supseteq A[\ell^n]\) for some \(n\geq 1\) then
we can find a polarization \(\lambda'\) on \(A\) such that
\(\lambda=\ell^n\lambda'\) by lemma \ref{lem:pol}(ii). Thus we
reduce to the case where \(\Lambda[\ell^\infty]\nsupseteq A[\ell^n]\)
for all \(n\geq 1\).

Now choose the minimal \(n\geq 1\) such that
\(\Lambda[\ell^\infty]\subseteq A[\ell^n]\), that is, such that
\(\Lambda[\ell^\infty]=\Lambda[\ell^n]\). Suppose that \(n>1\). Then
\(\ell\cdot\Lambda[\ell^2]\subseteq\Lambda[\ell]\), and
\(\ell\cdot\Lambda[\ell^2]\neq(0)\), for otherwise we would have
\(\ell^{n-2}\cdot\Lambda[\ell^n]\subseteq\Lambda[\ell^2]=\Lambda[\ell]\),
whence \(\Lambda[\ell^\infty]=\Lambda[\ell^n]=\Lambda[\ell^{n-1}]\),
contradicting the minimality of \(n\).

Let \(a,a'\in\ell\cdot\Lambda[\ell^2]\), and choose
\(b\in\Lambda[\ell^2]\) such that \(\ell b=a'\). Then
\(e^\lambda(a,a')=\ebar_\ell(a,\lambda(b))=1\), since
\(b\in\ker\lambda\). Hence, by lemma \ref{lem:pol}(iii), we
have an isogeny \(\pi\colon(A,\lambda)\to(B,\lambda')\) defined over
\(k\), with \(\ker\pi=\ell\cdot\Lambda[\ell^2]\). Then
\(\deg\lambda'\) is a proper divisor of \(\deg\lambda\) and, since the
action of \(\End_k(A)^\mathrm{s}\) preserves \(\ker\pi\), the isogeny
\(\pi\) induces \(\End_k(A)^\mathrm{s}\into\End_k(B)^\mathrm{s}\), by
lemma \ref{lem:desc}.

In this way we reduce to the case where
\(\Lambda[\ell^\infty]\) is properly contained in \(A[\ell]\). Since
the degree of a polarization is always a square, the result follows.
\end{Proof}

\bigskip
Finally, we note a result on the nature of those primes which divide
the degree of a polarization on an abelian surface with real
multiplication.

\begin{Lemma}\label{lem:sr}
Let \(R\) be an order in \(F\), and suppose that \((A,\lambda)\) is a
polarized abelian surface with an embedding
\(R\into\End(A)^\mathrm{s}\). Then every prime factor of
\(\deg\lambda\) which does not divide the conductor of \(R\) is either
split or ramified in \(F\).
\end{Lemma}

\begin{Proof}
Let \(H_{\Delta,d}\) denote the moduli space of triples
\((A,\lambda,i)\) formed of an abelian surface \(A\), a polarization
\(\lambda\) on \(A\) with \(\deg\lambda=d^2\), and an embedding
\(i\colon R\into\End(A)^\mathrm{s}\), where \(R\) is the order in
\(F\) of discriminant \(\Delta\). (Note that the discriminant
\(\Delta\) and conductor \(f(R)\) of \(R\) are linked by \(\Delta =
f(R)^2\mathrm{disc}(F)\).) In \cite[chapter IX, \S 2]{vdG}, van der
Geer describes \(H_{\Delta,d}\) as a quotient of the Siegel upper
half-space of degree 2, and shows that \(H_{\Delta,d}\neq\emptyset\)
precisely when \(\Delta\) is a square modulo \(4d\). The result as
stated is a straightforward deduction.
\end{Proof}

\section{Enlarging the endomorphism ring}
\label{sec:max}

In this section, we consider to what extent one may take isogenies to
enlarge the endomorphism ring of a given polarized abelian surface
\((A,\lambda)\) with real multiplication by a field \(F\) without
altering the degree of the polarization \(\lambda\). (It is well-known
that one may reduce to the case where \(\End(A)\) is the full ring of
integers of \(F\) if one disregards what happens to the polarization.)

\begin{Theorem}\label{thm:max}
Let \(\ell\) be an odd prime, and let \(R'\supset R\) be orders in
\(F\) with \([R':R]=\ell\). Further, let \((A,\lambda)\) be a
polarized abelian surface, defined over \(k\), with an embedding
\(i\colon R\into\End_k(A)^\mathrm{s}\), and suppose that
\(\ell\nmid\deg\lambda\).

Then there is an isogeny \(\pi\colon A\to B\) defined over \(k\) and a
polarization \(\lambda'\) on \(B\) also defined over \(k\) such that
\(\deg\lambda' = \deg\lambda\) and \(R'\into\End_k(B)^\mathrm{s}\).
\end{Theorem}

\begin{Proof}
We may clearly suppose that \(A\) has a symmetric action of \(R\), but
not of \(R'\).

Write \(R'=\Z[\alpha]\) for some \(\alpha\in F\), then
take \(\pi\colon A\to B\) to be an isogeny with kernel
\(\left(\ell\Z[\alpha]\right)\cdot A[\ell^{2}]\). (Note that we must
have an action of \((\ell\alpha)\) on \(A\).)

Now we have \(\ker\pi\supseteq A[\ell]\), and
\[
\ell\ker\pi = \left(\ell\Z[\alpha]\right)\cdot A[\ell] =
\im(\ell\alpha\mid A[\ell]).
\]
Thus \(\#\ker\pi = \ell^{4+t}\), where \(t =
\dim_{\F{\ell}}\im(\ell\alpha\mid A[\ell])\).

Set \(\lambda_A = \ell^{3}\lambda\). Then \(\ker\lambda_A\supseteq
A[\ell^{3}]\supseteq\ker\pi\). Also, if \(a,a'\in\ker\pi\) and
\(\ell^{2}b = a'\) then
\begin{eqnarray*}
e^{\lambda_A}(a,a')
\mytab = \mytab \ebar_{\ell^{2}}(a,\lambda_A(b)) \\
\mytab = \mytab \ebar_{\ell^{2}}(a,\ell\lambda(a')) \\
\mytab = \mytab \left[\ebar_{\ell^2}^\lambda(a,a')\right]^\ell \\
\mytab = \mytab \ebar_{\ell}^\lambda(\ell a, \ell a').
\end{eqnarray*}
But \(a = (\ell\eps)c\) and \(a'=(\ell\eps')c'\) for some
\(\eps,\eps'\in\Z[\alpha]\) and some \(c,c'\in A[\ell^{2}]\), so
\begin{eqnarray*}
e^{\lambda_A}(a,a')
\mytab = \mytab
\ebar_{\ell}^\lambda(\ell(\ell\eps)c,\ell(\ell\eps')c') \\
\mytab = \mytab
\ebar_{\ell}^\lambda(\ell^{2}(\ell\eps\eps')c, \ell c') \\
\mytab = \mytab 1.
\end{eqnarray*}

Hence by lemma \ref{lem:pol}(iii), we may write
\(\lambda_A = \pi^{*}\lambda'\) for some polarization \(\lambda'\) on
\(B\), defined over \(k\). Note that
\[
\deg\lambda' = \frac{\ell^{12}\deg\lambda}{\ell^{2(4+t)}}
             = \ell^{2(2-t)}\deg\lambda.
\]
Further, applying lemma \ref{lem:desc} to
\(\pi\colon(A,\lambda_A)\to(B,\lambda')\) shows that \(B\) has an action of
\(\Z[\ell\alpha]\) which is defined over \(k\) and fixed by the
Rosati involution.

We can define an action of \(\alpha\) on \(B\) as follows: given
\(b\in B(\kbar)\), choose \(a\in A(\kbar)\) such that \(\ell\pi(a)=b\)
and set \(\alpha\cdot b = \pi\big((\ell\alpha)\cdot a\big)\). To check
that this action is well-defined, we must check that
\(a\in\ker(\ell\pi)\) implies \((\ell\alpha)\cdot a\in\ker\pi\). But if
\(a\in\ker(\ell\pi)\), then
\[
\ell(\ell\alpha)\cdot a = (\ell\alpha)\cdot(\ell a) \in
(\ell\alpha)\cdot\ker\pi \subseteq \ell\cdot\ker\pi,
\]
by the definition of \(\ker\pi\). Hence \((\ell\alpha)\cdot a\in
(\ker\pi + A[\ell]) = \ker\pi\).

One now readily checks that this action is defined over \(k\) and
fixed by the Rosati involution. All that remains, then, is to show
that \(\deg\lambda' = \deg\lambda\), that is, \(t=2\).

Note that \(e_\ell^\lambda\) is a nondegenerate pairing on \(T_\ell A\)
since \(\ell\nmid\deg\lambda\). Also, \((\ell\alpha)\) is not zero
on \(A[\ell]\) since we are assuming that \(A\) does not have an
action of \(\Z[\alpha]\). However, \((\ell\alpha)^2\) does act as
zero on \(A[\ell]\) and hence, applying lemma \ref{lem:easy}(ii), we
conclude that \(t = 2\) as required.
\end{Proof}

\bigskip
We note that the proof of theorem \ref{thm:max} does not generalize to
higher dimensional varieties. Even if one generalizes the proof to
accommodate the fact that the orders in a number field \(F\) of degree
greater than 2 may not be of the form \(\Z[\alpha]\), there remains a
further obstacle to the method. This is that not only can one no
longer guarantee that the degree of the new polarization is the same
as that of the old, but when the variety has odd dimension then the
degree of the new polarization must be greater than that of the
old. If one follows through the proof with \(\dim A = g\) (rather than
\(\dim A=2\)), then one produces a polarization \(\lambda'\) on the
variety \(B\) such that
\[
\deg\lambda' = \ell^{2(g-t)}\deg\lambda,
\]
with \(t\) defined as before. We can argue, as before, that \(t\neq
0\), but we know from lemma \ref{lem:easy} that \(t\) must be
even. Hence we must have \(\deg\lambda' > \deg\lambda\) when \(g\) is
odd.

\section{Reducing the degree of a polarization}
\label{sec:sr}

In this final section, we consider to what extent one might be able to
reduce the degree of a given polarization by taking isogenies. The
result below can be argued at least heuristically from the isomorphism
between \(\mathrm{NS}(A)\otimes\Q\) and the algebra of symmetric
elements of \(\End(A)\otimes\Q\), but we wish to take care that we can
find an isogenous variety \(B\), with a polarization
\(\lambda'\) on \(B\) (not just a homogeneous polarization) such that
\(\deg\lambda'\) divides \(\deg\lambda\), and that
\(\End_k(A)^\mathrm{s}\) embeds into \(\End_k(B)^\mathrm{s}\).

\begin{Theorem}\label{thm:sr}
Let \(R\) be an order in \(F\), and let \((A,\lambda)\) be a polarized
abelian surface over \(k\), with an embedding \(i\colon
R\into\End_k(A)^\mathrm{s}\). Suppose that \(\ell\) is an odd prime
factor of \(\deg(\lambda)\) which does not divide the conductor of
\(R\), and which is reducible in \(R\).

Then there is an isogeny \(\pi\colon A\to B\) defined over \(k\) and a
polarization \(\lambda'\) on \(B\), defined over \(k\), such that
\(\deg\lambda'\) divides \(\ell^{-2}\deg\lambda\), and such that \(R\)
embeds into \(\End_k(B)^\mathrm{s}\).
\end{Theorem}

\begin{Proof}
Write \(\Lambda=\ker\lambda\). After lemma \ref{lem:sqf} we are
reduced to considering the case where \(\Lambda[\ell^\infty]\subset
A[\ell]\) and \(\dim_{\F{\ell}}(\Lambda[\ell^\infty]) = 2\).

\bigskip
We are given that \(\ell = \alpha_1\alpha_2\) in \(R\), where
\(|\norm{F/\Q}{\alpha_1}| = |\norm{F/\Q}{\alpha_2}| = \ell\).
The elements \(\alpha_1\) and \(\alpha_2\) may or may not be
associates in the ring of integers \(\mO\) of \(F\); we consider these
cases separately.

If \(\alpha_1\) and \(\alpha_2\) are not associates in \(\mO\), then
the ideals \(\alpha_1\mO\) and \(\alpha_2\mO\) are coprime, and hence
we may write
$$
f = \alpha_1\beta_1 + \alpha_2\beta_2 \eqno{(*)}
$$
in \(R\), where \(f\) is the conductor of \(R\). By hypothesis, \(f\)
is invertible modulo \(\ell\). Hence equation \((*)\)
induces decompositions \(A[\ell]=A[\alpha_1]\oplus A[\alpha_2]\) and
\(\Lambda[\ell]=\Lambda_1\oplus\Lambda_2\), where \(\Lambda_i =
\Lambda\cap A[\alpha_i]\) for \(i=1,2\). Note that
\(\dim_{\F{\ell}}A[\alpha_1]= 2 = \dim_{\F{\ell}}A[\alpha_2]\).

Suppose that both \(\Lambda_1\) and \(\Lambda_2\) are nonzero. Then
\(\Lambda_1\) is 1-dimensional and rationally defined. The pairing
\(e^\lambda\) is alternating, so will be trivial on
\(\Lambda_1\times\Lambda_1\) and we may take an isogeny
\(\pi\colon(A,\lambda)\to(B,\lambda')\) with \(\ker\pi=\Lambda_1\) by
lemma \ref{lem:pol}(iii). Note that \(\deg\lambda' =
\ell^{-2}\deg\lambda\). Also, \(R\) preserves \(\Lambda_1\), and so
\(R\into\End_k(B)^\mathrm{s}\) by lemma \ref{lem:desc}.

Suppose instead that one of the factors is zero. Without loss of
generality, then, we have \(\Lambda[\ell]=\Lambda_1\). But, simply
considering group orders, this forces \(\Lambda_1 =
A[\alpha_1]\). Hence, we may apply lemma \ref{lem:pol}(ii) to
give a polarization \(\lambda'\) on \(A\) such that
\(\lambda=\lambda'\alpha_1\). Again,
\(\deg\lambda'=\ell^{-2}\deg\lambda\).

\bigskip
Now suppose that \(\alpha_1\) and \(\alpha_2\) are associates in
\(\mO\). Then, on writing \(\Lambda_1=\Lambda\cap A[\alpha_1]\) as
before, we have \(\Lambda[\ell]\supseteq\Lambda_1\supseteq(0)\). Note
that now \(\alpha_1\cdot\Lambda[\ell]\subseteq\Lambda_1\) and so
\(\Lambda_1\neq(0)\), since \(\Lambda[\ell]\neq(0)\) (arguing
similarly to the third paragraph of the proof of lemma
\ref{lem:sqf}). We now proceed exactly as in the previous case, since
\(\Lambda\cap A[\alpha]\) is of dimension 1 or 2 over \(\F{\ell}\).
\end{Proof}

\bigskip
There is an analogous statement to the theorem above for an abelian
variety of general dimension \(g\): the same method allows one to remove
from \(\deg\lambda\) those prime factors \(\ell\) which split in \(R\)
as a product of \(g\) (not necessarily distinct) elements of
\(R\) each of norm \(\pm\ell\). Lemma \ref{lem:sr} guarantees that
this is a sensible condition to impose on the prime \(\ell\) when
dealing with abelian surfaces, but it is not clear to the author
whether or not lemma \ref{lem:sr} generalizes in the appropriate
fashion to varieties of dimension \(3\) or more.

\bigskip
We note that in the statement of theorem \ref{thm:sr} we require not
just that the prime factors of \(\deg\lambda\) are no longer prime in
\(R\), but that they should be reducible in \(R\) (which is, of
course, a strictly stronger condition in general). This is why, in
theorem \ref{thm:main}, we require \(F\) to have class number 1. We
now present an argument to show why we might expect that this
condition is necessary in general.  (One should note that, in some
sense, the case where \(F\) has class number greater than 1 seems to
be a rarity: Wang \cite{Wang} notes that no simple 2-dimensional
factor of \(J_0(N)\) with \(N\leq 500\) has real multiplication by a
field with class number greater than 1.)

Consider the following hypothetical example. Suppose that
\((A,\lambda)\) is a polarized abelian surface over \(k\), with
\(\End_k(A)^\mathrm{s}\otimes\Q \cong F\), a real quadratic field
(perhaps \(k=\Q\) and \(A\) is a simple factor of \(J_0(N)\) for some
\(N\)). Suppose also that \(\ell\) is a prime such that
\(\ell^2 \parallel\deg\lambda\), and that we have a \(k\)-isogeny
\(\pi\colon A\to B\), and a polarization \(\lambda'\) on \(B\) such
that \(\deg\lambda'\) divides \(\ell^{-2}\deg\lambda\). We may set
\(\lambda_A := \pi^{*}\lambda'\), and then \(\lambda =
\lambda_A\circ\alpha\) for some \(\alpha\in F\). Note that this forces
\[
\deg\lambda = \left(\norm{F/\Q}{\alpha}\right)^2
\left(\deg\pi\right)^2 \left(\deg\lambda'\right),
\]
which implies that
\[
\left(\norm{F/\Q}{\alpha}\right)\left(\deg\pi\right) =
\ell m, \textrm{ where } m^2 = \frac{\deg\lambda}{\ell^2\deg\lambda'};
\]
thus one or other of \(\deg\pi\) and \(\norm{F/\Q}{\alpha}\) must be
divisible by an odd power of \(\ell\). 

The degree of \(\pi\) can be divisible by an odd power of \(\ell\)
only if there is a suitable \(k\)-rational subgroup of \(A\)
contained in \(A[\ell^\infty]\) but not equal to any \(A[\ell^r]\),
which need not happen in general.

On the other hand, if \(\deg\pi\) is divisible by an even power of
\(\ell\), then \(\norm{F/\Q}{\alpha}\) must be divisible
by an odd power of \(\ell\), and every other prime factor of
\(\norm{F/\Q}{\alpha}\) which appears to a positive power must also be
a factor of \(\deg\lambda\). In the extreme case, when
\(\deg\lambda=\ell^2\), this forces \(\norm{F/\Q}{\alpha}=\ell\),
which need not be possible if \(F\) has class number greater than 1.

\bigskip
We conclude with a deduction of theorem \ref{thm:main} (which we re-state here
for convenience) from theorems \ref{thm:max} and \ref{thm:sr}.

\setcounter{Theorem}{0}
\begin{Theorem}
Let \(R\) be an order in \(F\) and let \((A,\lambda)\) be a polarized
abelian surface over \(k\) with an embedding \(i\colon
R\into\End_k(A)^\mathrm{s}\). Suppose that \(F\) has class number
\(1\), and that the degree of \(\lambda\) and the conductor of \(R\)
are odd and coprime.

Then there is a \(k\)-isogenous abelian surface \(B\) with a principal
polarization \(\lambda'\) on \(B\) and a symmetric action of the ring
of integers of \(F\), both defined over \(k\).
\end{Theorem}

\begin{Proof}
Theorem \ref{thm:max} applied repeatedly to the prime factors of
the conductor of \(R\) guarantees the existence of a \(k\)-isogenous
abelian surface \(A_1\) with a polarization \(\lambda_1\) on \(A_1\)
satisfying \(\deg\lambda_1=\deg\lambda\), such that
\(\End_k(A_1)^\mathrm{s}\) contains a copy of the ring of integers of
\(F\).

Applying lemma \ref{lem:sr} and the hypothesis that \(F\) has class
number 1, every prime factor of \(\deg\lambda_1\) is reducible in the
ring of integers of \(F\). Thus theorem \ref{thm:sr} applied
repeatedly to the prime factors of \(\deg\lambda_1\) guarantees the
existence of an abelian surface \(B\) and polarization \(\lambda'\) on
\(B\) as required.
\end{Proof}

\begin{flushright}
\textit{John Wilson, Magdalen College, Oxford OX1 4AU}
\par\texttt{wilsonj@maths.ox.ac.uk}
\end{flushright}

\end{document}